\documentstyle[11pt]{article}

\title{THE GALAXIES OF NONSTANDARD ENLARGEMENTS OF TRANSFINITE GRAPHS
OF HIGHER RANKS: II}
\author{A. H. Zemanian}
\date{}

\begin{document}
\newcommand{\N} {I \kern -4.5pt N}
\newcommand{\R} {I \kern -4.5pt R}
\newcommand{\Z} {Z \kern -7.5pt Z}
\maketitle
\baselineskip21pt

{\ Abstract --- This report is an improvement of a prior report
\cite[Report 814]{814}.
It sharpens the principal theorem (Theorem 5.1 of Report 814)
and also simplifies its proof.  There are also several minor changes 
involving clarifications and corrections of misprints.
The prior abstract remains the same as follows:

In a previous work, the galaxies of the nonstandard enlargements 
of connected, conventionally infinite graphs and also 
of connected transfinite graphs 
of the first rank of transfiniteness were 
defined, examined, and illustrated by some examples.  In this work it is 
shown how the results of the prior 
work extend to transfinite graphs 
of higher ranks. Among those results are following:
Any such enlargement either has exactly one galaxy, its principal one, or 
it has infinitely many such galaxies.  In the latter case, the galaxies are 
partially ordered by there ``closeness'' to the principal galaxy.
Also, certain sequences of galaxies whose members are totally ordered by 
that ``closeness'' criterion are identified.

} 

\section{Introduction}

In some previous works, the ideas of ``nonstandard graphs''
\cite[Chapter 8]{gn} and the ``galaxies of nonstandard enlargements
of graphs'' \cite{gal}, \cite{819} were defined and examined.  However, 
all this was done only for conventionally infinite graphs and
transfinite graphs of the first rank of transfiniteness.
The purpose of this work is to define and examine nonstandard
transfinite graphs of higher ranks of transfiniteness.

This paper is written as a sequel to \cite{819} and uses a 
symbolism and terminology consistent with that prior work.
We also use a variety of results concerning 
transfinite graphs, and these may all be found in \cite{gn}. 
We refer the reader to those sources for such information.
For instance, the ``hyperordinals'' are constructed in 
much the same way as are the hypernaturals, and their definitions are given 
in \cite[Section 5]{819}.  Furthermore, we use herein the 
idea of ``wgraphs,'' which are transfinite graphs based upon walks
rather than paths.  This avoids some of the difficulties associated with
path-based transfinite graphs and is in fact both simpler 
and more general than the path-based theory of transfinite graphs.
Wgraphs and their transfinite extremities are defined and examined in 
\cite[Sections 5.1 to 5.6]{gn}.  Lengths of walks on wgraphs 
and ``wdistances'' based on such lengths are discussed 
in \cite[Sections 5.7 and 5.8]{gn}.
All this is assumed herein as being known.
It is also assumed that all the wgraphs considered herein are
wconnected.

Our arguments will be based on ultrapower constructions, and, to this end, 
we assume throughout that a free ultrafilter $\cal F$ has been 
chosen and fixed.  Finally, when adding ordinals, we always take it that 
ordinals are in normal form and that the natural summation of ordinals is 
being used \cite[pages 354-355]{ab}.

It is a fact about a 
transfinite wgraph $G^{\nu}$ of rank 
$\nu$ that it contains wsubgraphs of all ranks 
$\rho$ with $0\leq\rho\leq\nu$, called $\rho$-wsections,
that at each rank $\rho$ the $\rho$-wsections are $\rho$-wgraphs by 
themselves and induce a partitioning of the branch set of 
$G^{\nu}$, and that the one and only
$\nu$-wsection is $G^{\nu}$ itself.
We define the ``enlargements'' $^{*}\!G^{\nu}$ of a transfinite 
wgraph $G^{\nu}$ and of its $\rho$-wsections 
in the next section. The galaxies of all ranks in $^{*}\!G^{\nu}$
are defined in Section 3.  A galaxy of rank $\rho$ ($0\leq \rho\leq\nu$)
is called a ``$\rho$-galaxy.''  

Within the enlargement 
$^{*}\!S^{\rho}$ of an $\rho$-wsection $S^{\rho}$ of $G^{\nu}$,
there is either 
exactly one $\rho$-galaxy, the ``principal $\rho$-galaxy,''
or infinitely many $\rho$-galaxies in addition to the principal $\rho$-galaxy.
The latter case arises 
when $G^{\nu}$ is locally finite in a certain way (Section 4),
but it may arise in other ways as well.  Moreover, the enlargements of 
all the $\rho$-wsections within a $(\rho +1)$-wsection
lie within the principal $(\rho+1)$-galaxy
of the enlargement of that $(\rho+1)$-wsection,
and so on through the wsections of higher ranks.
In that latter case still, there
will be a two-way infinite sequence of $\rho$-galaxies that are totally 
ordered according to their ``closeness to the principal $\rho$-galaxy,''
and there may be many such totally ordered sequences of 
$\rho$-galaxies.  When there are many $\rho$-galaxies in $^{*}\!S^{\rho}$,
they are partially ordered, again according to their closeness
to the principal $\rho$-galaxy (Section 5).

\section{Enlargements of $\nu$-Graphs and Hyperdistances in the Enlargements}

First of all, the enlargement $^{*}\!G^{0}=\{\,^{*}\!X^{0},\,^{*}\!B\,\}$
of a conventionally infinite 0-graph and the enlargement 
$^{*}\!G^{1}=\{\,^{*}\!X^{0},\,^{*}\!B,\,^{*}\!X^{1}\,\}$
of a transfinite 1-graph are discussed in \cite[Sections 2 and 8]{gal}.
These prior constructs will be encompassed by the more
general development we now undertake.  We shall assume that the
rank $\nu$ is no larger than $\omega$.  The extensions to higher ranks
of transfiniteness proceeds in much the same way.

Consider a wconnected transfinite wgraph of rank $\nu$ ($0\leq\nu\leq\omega$):
\[ G^{\nu}\;=\;\{X^{0},B,X^{1},\ldots,X^{\nu}\} \]
where $X^{0}$ is a set of 0-nodes, $B$ is a set of branches
(i.e., two-element sets of 0-nodes), and $X^{\rho}$
$(\rho=1,\ldots,\nu)$ is a set of $\rho$-wnodes.  It is assumed 
that each $X^{\rho}$ $(\rho=0,\ldots,\nu)$ is nonempty except possibly for 
$\rho=\vec{\omega}$.  In general, $X^{\vec{\omega}}$ may be empty.

The ``enlargement'' $^{*}\!G^{\nu}$ of $G^{\nu}$ is defined as follows: 
Two sequences $\langle x_{n}^{\rho}\rangle$ and 
$\langle y_{n}^{\rho}\rangle$ of $\rho$-wnodes in $G^{\nu}$ 
(i.e., $x_{n}^{\rho},y_{n}^{\rho}\in X^{\rho}$)
are taken to be equivalent if $\{ n\!: x_{n}^{\rho}=y_{n}^{\rho}\}\in{\cal F}$.
This is truly an equivalence 
relation on the set of sequences of $\rho$-wnodes, 
as is easily shown.  Each equivalence class ${\bf x}^{\rho}$ will
be called a $\rho$-{\em hypernode} and will be represented by 
${\bf x}^{\rho}=[x_{n}^{\rho}]$ where the $x_{n}^{\rho}$ 
are elements of any one 
of the sequences in that equivalence class.  
We let $^{*}\!X^{\rho}$ denote the set of all such equivalence classes
(i.e., the set of all $\rho$-hypernodes).
Then, the {\em enlargement} $^{*}\!G^{\nu}$ of $G^{\nu}$ is 
the set 
\[ ^{*}\!G^{\nu}\;=\;\{\,^{*}\!X^{0},\,^{*}\!B,\,^{*}\!X^{1}\ldots,
\,^{*}\!X^{\nu}\,\}. \]
The elements of $\,^{*}B$ are called {\em hyperbranches} and have 
been defined in \cite[Section 8.1]{gn}.  Here, too, 
$\,^{*}\!X^{\rho}$ is nonempty if $\rho\neq \vec{\omega}$; 
$\,^{*}\!X^{\vec{\omega}}$ may be empty.

Next, we wish to define the ``hyperdistances''  between the hypernodes of 
$^{*}\!G^{\nu}$.  The ``length'' $|W_{xy}|$ of any 
two-ended walk $W_{xy}$ terminating at two wnodes
$x$ and $y$ of any ranks in $G^{\nu}$ is defined in \cite[Section 5.7]{gn}.
Also, the wdistance $d(x,y)$ between those two wnodes is
\[ d(x,y)\;=\;\min\{|W_{x,y}|\} \]
where the minimum is taken over the lengths $|W_{xy}|$ of all
the walks $W_{xy}$ in $G^{\nu}$ that terminate at $x$ and $y$.
The minimum exists because those lengths comprise a 
well-ordered set of ordinals.  Furthermore, 
a wnode is said to be {\em maximal} if it is not embraced by a 
wnode of higher rank.  In these definitions, $x$ and $y$ may 
be either maximal or nonmaximal wnodes.  Note that the wdistance 
measured from any nonmaximal wnode $z$ is the same as that measured from the 
maximal wnode $x$ that embraces $z$.  We also set $d(x,x)=0$.
Thus, $d$ is an ordinal-valued metric defined on the maximal
wnodes in $\cup_{\rho=0}^{\nu} X^{\rho}$.  (The axioms of a metric
are readily verified.)

Given two hypernodes ${\bf x}=[x_{n}]$ and ${\bf y}=[y_{n}]$ of any ranks
in $^{*}\!G^{\nu}$, we defined the {\em hyperdistance}
${\bf d}$ between them as the internal function
\[ {\bf d}({\bf x},{\bf y})\;=\;[ d(x_{n},y_{n}) ]. \]
We say that a hypernode ${\bf x}^{\rho}=[x_{n}^{\rho}]$ is {\em maximal}
if it is not embraced by a hypernode ${\bf y}^{\gamma}=[y_{n}^{\gamma}]$
of higher rank $(\gamma > \rho)$ (i.e., $x_{n}^{\rho}$ 
is not embraced by $y_{n}^{\gamma}$ for almost all $n$).
Upon restricting ${\bf d}$ to 
the maximal hypernodes in $^{*}\!G^{\nu}$, we have that
this restricted ${\bf d}$ satisfies the metric axioms except that it 
is hyperordinal-valued.  In particular, we have by the transfer 
principle that the triangle inequality holds for any three maximal
hypernodes ${\bf x}$, ${\bf y}$, and ${\bf z}$, namely,
\begin{equation}
{\bf d}({\bf x},{\bf z})\;\leq\; {\bf d}({\bf x},{\bf y})\:+\: {\bf d}({\bf y},{\bf z}).  \label{2.1}
\end{equation}

\section{The Galaxies of $^{*}\!G^{\nu}$}

We continue to assume that the rank $\nu$ is no larger than $\omega$.
Also, we assume at first that the rank $\rho$ is a natural number no larger
than $\nu$.  Consider the $\nu$-wgraph $G^{\nu}$ and its enlargement 
$^{*}\!G^{\nu}$.  Two hypernodes ${\bf x}=[x_{n}]$ and 
${\bf y}=[y_{n}]$ of any ranks in 
$^{*}\!G^{\nu}$ will be said to be in the same 
{\em nodal $\rho$-galaxy} $\dot{\Gamma}^{\rho}$ if there exists
a natural number $\mu_{xy}$ depending on ${\bf x}$ and ${\bf y}$ such that
$\{n\!: d(x_{n},y_{n})\leq \omega^{\rho}\cdot \mu_{xy}\}\in {\cal F}$.
In this case, we say that ${\bf x}$ and ${\bf y}$ 
are $\rho$-{\em limitedly distant}.
This defines an equivalence relation on the set 
$\cup_{\gamma=0}^{\nu}$$^{*}\!X^{\gamma}$ of all the hypernodes in 
$^{*}\!G^{\nu}$, and thus $\cup_{\gamma=0}^{\nu}$$^{*}\!X^{\gamma}$ is
partitioned into nodal $\rho$-galaxies.
The proof of this is the same as that given in \cite[Section 9]{gal}
except that the rank 1 therein is now replaced by $\rho$.
By the same arguments, we have that, for $\alpha\leq \rho$, the 
nodal $\alpha$-galaxies provide a finer partitioning of 
$\cup_{\gamma=0}^{\nu}$$^{*}\!X^{\gamma}$ than do the nodal $\gamma$-galaxies.
Moreover, any nodal $\rho$-galaxy is partitioned by the 
nodal $\alpha$-galaxies $(\alpha < \rho)$ in that nodal $\rho$-galaxy.

Corresponding to each nodal $\rho$-galaxy $\dot{\Gamma}^{\rho}$, 
we define a $\rho$-{\em galaxy} $\Gamma^{\rho}$ of 
$^{*}\!G^{\nu}$ as the nonstandard wsubgraph of $^{*}\!G^{\nu}$
induced by all the $\rho$-hypernodes in $\dot{\Gamma}^{\rho}$;
that is, along with the hypernodes of $\dot{\Gamma}^{\rho}$,
we have hyperbranches whose incident 0-hypernodes are in $\dot{\Gamma}^{\rho}$.
Note that every hyperbranch must lie in a single $\rho$-galaxy
for every $\rho$ because their incident 0-hypernodes are at a 
hyperdistance of 1.

Let us now turn to the case where $\rho=\vec{\omega}$.  Now, $\nu$ is 
either $\vec{\omega}$ or $\omega$.  The definition of 
the $\vec{\omega}$-galaxies
is rather different.  Two hypernodes ${\bf x}=[x_{n}]$
and ${\bf y}=[y_{n}]$ in $^{*}\!G^{\nu}$ will be said to be in the
same $\vec{\omega}$-{\em galaxy} $\dot{\Gamma}^{\vec{\omega}}$ 
if there exists a natural number
$\mu_{xy}$ depending on ${\bf x}$ and ${\bf y}$ such that 
$\{n\!: d(x_{n},y_{n})\leq \omega^{\mu_{xy}}\}\in {\cal F}$.
Thus, ${\bf x}$ and ${\bf y}$ are in $\dot{\Gamma}^{\vec{\omega}}$
if and only if $\{n\!: d(x_{n},y_{n}) < \omega^{\omega}\}\in {\cal F}$.
In this case, we say that 
${\bf x}$ and ${\bf y}$ are $\vec{\omega}$-{\em limitedly distant}.
Here, too, the property of being $\vec{\omega}$-limitedly distant
defines an equivalence relation on the set of all 
hypernodes in $^{*}\!G^{\nu}$.  So, the nodal $\vec{\omega}$-galaxies 
partition the set of all hypernodes.  Then, the $\vec{\omega}$-galaxy
$\Gamma^{\vec{\omega}}$ corresponding to any 
nodal $\vec{\omega}$-galaxy $\dot{\Gamma}^{\vec{\omega}}$
consists of the hypernodes in $\dot{\Gamma}^{\vec{\omega}}$ along with 
the hyperbranches whose 0-hypernodes are in $\dot{\Gamma}^{\vec{\omega}}$.

Finally, the $\omega$-galaxies of a nonstandard $\omega$-wgraph
$^{*}\!G^{\omega}$ are defined just as are the $\rho$-galaxies 
of natural-number ranks.  We now require that $\{n\!: d(x_{n},y_{n}\leq
\omega^{\omega}\cdot \mu_{xy}\}\in {\cal F}$ for some natural number
$\mu_{xy}$ depending upon ${\bf x}=[x_{n}]$ and ${\bf y}=[y_{n}]$
in order for ${\bf x}$ and ${\bf y}$ to be in the same nodal $\omega$-galaxy.
When this is so, we again say that ${\bf x}$ and ${\bf y}$ are
$\omega$-{\em limitedly distant}.  The same partitioning properties hold.

In general now, let $G^{\nu}$ be a $\nu$-wgraph where
$0\leq\nu\leq\omega$, possibly $\nu=\vec{\omega}$.
The {\em principal $\nu$-galaxy} $\Gamma_{0}^{\nu}$ of
$^{*}\!G^{\nu}$ is that $\nu$-galaxy whose hypernodes are $\nu$-limitedly 
distant from a standard hypernode of $^{*}\!G^{\nu}$.  
We shall show later on that $^{*}\!G^{\nu}$ either has exactly one 
$\nu$-galaxy, its principal one,
or has infinitely many of them.   

Let us now recall another 
definition concerning standard transfinite wgraphs.  
A $\rho$-{\em wsection} $S^{\rho}$ of $G^{\nu}$ $(\rho <\nu)$
is a maximal $\rho$-wsubgraph of the $\rho$-wgraph of
$G^{\nu}$ that is $\rho$-wconnected.
This {\em $\rho$-wconnectedness} means that, for every two wnodes in 
$S^{\rho}$, there is a two-ended $\alpha$-walk 
$(\alpha\leq\rho)$ terminating at those wnodes.  (When $\rho=\vec{\omega}$, 
we have that $\alpha <\vec{\omega}$.) 
Furthermore, the branch set of any $\rho$-wsection
$S^{\rho}$ is partitioned by the branch sets of the $\alpha$-wsections
lying in $S^{\rho}$.

Consider now any $(\rho-1)$-wsection $S^{\rho-1}$ within a 
$\rho$-wsection $S^{\rho}$.  A $\rho$-wnode (not necessarily maximal now)
will be called {\em incident} to $S^{\rho-1}$ if $x^{\rho}$ embraces 
an $\alpha$-wtip $t^{\alpha}$ $(\alpha<\rho)$ belonging to $S^{\rho-1}$
(that is, $S^{\rho-1}$ contains a one-ended extended $\alpha$-walk 
that is a representative of $t^{\alpha}$, 
where an ``extended'' $\alpha$-walk is defined 
in \cite[page 74]{gn}).  This definition includes 
the case where the $\alpha$-walk is simply a single branch;  we let
the extremity of the branch be an $(-1)$-{\em tip}.  Also, when 
$\rho=\omega$, $\rho-1$ denotes $\vec{\omega}$.

{\bf Lemma 3.1.}  {\em Given any two $\rho$-wnodes $x^{\rho}$ and $y^{\rho}$
incident to a $(\rho-1)$-wsection, there exists an $\alpha$-walk
$(\alpha\leq\rho-1)$ in $S^{\rho-1}$ that reaches $x^{\rho}$ and $y^{\rho}$.}

{\bf Proof.}  Let $W_{x}^{\alpha_{1}}$ (resp. $W_{y}^{\alpha_{2}}$)
be a representative walk for the $\alpha_{1}$-wtip
(resp. $\alpha_{2}$-wtip) embraced by
$x^{\rho}$ (resp. $y^{\rho}$). We have that 
$\alpha_{1}, \alpha_{2}\,\leq\,\rho-1$. 
Let $u$ (resp. $v$) be a wnode of $W_{x}^{\alpha_{1}}$ 
(resp. $W_{y}^{\alpha_{2}}$) not embraced by $x^{\rho}$ (resp. $y^{\rho}$).
By the $(\rho-1)$-wconnectedness of $S^{\rho-1}$, there is a two-ended walk
$W_{uv}$ of rank no larger than $\rho-1$ that terminates at $u$ and $v$.
Then, the walk that passes first from $x^{\rho}$ along $W_{x}^{\alpha_{1}}$
to $u$, then along $W_{uv}$ to $v$, and finally along 
$W_{y}^{\alpha_{2}}$ to $y^{\rho}$ is the asserted walk.  $\Box$

Now, each $\rho$-wsection $S^{\rho}$ of $G^{\nu}$ is a $\rho$-wgraph by
itself, and therefore the enlargement $^{*}\!S^{\rho}$ of $S^{\rho}$
has its own principal $\rho$-galaxy $\Gamma_{0}^{\rho}(S^{\rho})$.

{\bf Theorem 3.2.}  {\em If $\alpha<\rho\leq \nu$ and if $S^{\alpha}$ 
is an $\alpha$-wsection 
lying $S^{\rho}$, then the enlargement $^{*}\!S^{\alpha}$ of $S^{\alpha}$
lies within the principal $\rho$-galaxy $\Gamma_{0}^{\rho}(S^{\rho})$
of $^{*}\!S^{\rho}$.}

{\bf Note.}  The conclusion means that every hypernode in 
$^{*}\!S^{\alpha}$ is a hypernode in $\Gamma_{0}^{\rho}(S^{\rho})$,
and consequently every hyperbranch in $^{*}\!S^{\alpha}$
is a hyperbranch in $\Gamma_{0}^{\rho}(S^{\rho})$.

{\bf Proof.}  Let ${\bf x}=[\langle x,x,x,\ldots\rangle]$
be any standard hypernode in $\Gamma_{0}^{\rho}(S^{\rho})$,
and let ${\bf y}=[\langle y,y,y,\ldots\rangle]$ be any 
standard hypernode in the  
principal $\alpha$-galaxy $\Gamma_{0}^{\alpha}(S^{\alpha})$
of $^{*}\!S^{\alpha}$.  Since $S^{\alpha}$ lies in $S^{\rho}$,
the standard wnodes $x$ and $y$ corresponding to ${\bf x}$ and ${\bf y}$
are no further apart than the wdistance $\omega^{\rho}\cdot k$
for some $k\in\N$.  (Indeed, there is a walk wconnecting them that does not
pass through any wnode of rank greater than $\rho$.)  Consequently,
${\bf x}$ and ${\bf y}$ are $\rho$-limitedly distant.  
Also, for every hypernode 
${\bf z}=[z_{n}]$ in $^{*}\!S^{\alpha}$, ${\bf y}$ and ${\bf z}$
are $\alpha$-limitedly distant, which implies that they are $\rho$-limitedly
distant.  So, by the triangle inequality (\ref{2.1}), ${\bf x}$ and ${\bf z}$
are $\rho$-limitedly distant.  Thus, ${\bf z}$ is in 
$\Gamma_{0}^{\rho}(S^{\rho})$.  $\Box$

Note that $^{*}\!S^{\alpha}$ may (but need not) have other $\alpha$-galaxies
besides its principal one $\Gamma_{0}^{\alpha}(S^{\alpha})$, and 
$^{*}\!S^{\rho}$ may have still other $\alpha$-galaxies not in
$\Gamma_{0}^{\rho}(S^{\rho})$.  Also, there may be a 
$\rho$-galaxy (possibly many of them) consisting of a single 
$\lambda$-hypernode when $\lambda >\rho$.  

Furthermore, it is possible of $^{*}\!G^{\nu}$ to have  exactly one
$\rho$-galaxy.  This occurs, for  instance, when there is a 
single node in $G^{\nu}$ to which all the other nodes of $G^{\nu}$
are connected through two-ended $\rho$-paths, with each such path
being connected to the rest of $G^{\mu}$ only at its terminal nodes.
When $^{*}\!G^{\nu}$
has exactly one $\rho$-galaxy, then $^{*}\!G^{\nu}$ has exactly one 
$\sigma$-galaxy for every $\sigma$ such that $\rho<\sigma\leq \nu$
because, if two hypernodes are $\rho$-limitedly distant, then they are 
also $\sigma$-limitedly distant. 

In view of all this, we again observe that the galaxies of 
$^{*}\!G^{\nu}$ can have rather complicated structures and dissimilarities.

\section{Locally Finite Sections and a Property of Their Enlargements}

In this and the next section, the rank $\rho$ is not allowed to be 
$\vec{\omega}$.  We now establish a sufficient (but not necessary) condition 
under which the enlargement $^{*}\!S^{\rho}$ has at least one 
$\rho$-galaxy different from its principal galaxy
$\Gamma_{0}^{\rho}(S^{\rho})$.  Let us first recall some definitions
for standard wgraphs.  

Assume initially that $\rho$ is a natural number.  Two 
$\rho$-wnodes of $S^{\rho}$ will be called $\rho$-{\em wadjacent} 
if they are incident to the same $(\rho-1)$-wsection.
A $\rho$-wnode will be called a {\em boundary $\rho$-wnode} if it is 
incident to two or more $(\rho-1)$-wsections.
A $\rho$-wsection $S^{\rho}$ will be called {\em locally 
$\rho$-finite} if each of its $(\rho-1)$-wsections has only
finitely many incident boundary $\rho$-wnodes.  
These same definitions hold when $\rho=\omega$ except that $\rho-1$ is 
understood to be $\vec{\omega}$.  
The case where $\rho=\vec{\omega}$ is prohibited in
the statements of this section.

In the following, we let $\rho$ be a natural number or 
$\rho=\omega$. When $\rho=\omega$, $\rho-1$ denotes $\vec{\omega}$.

{\bf Lemma 4.1.}  {\em Let $x^{\rho}$ be a boundary $\rho$-wnode.
Then, any $\rho$-walk that passes through $x^{\rho}$ from one 
$(\rho-1)$-wsection $S_{1}^{\rho-1}$ incident to $x^{\rho}$
to another $(\rho-1)$-wsection 
$S_{2}^{\rho-1}$ incident to $x^{\rho}$ must have a length no less
that $\omega^{\rho}$ (thus, when $\rho=\omega$, a length no less than 
$\omega^{\omega}$).}

{\bf Proof.}  
The only way such a walk can have a length less 
than $\omega^{\rho}$ is if it avoids traversing a 
$(\rho-1)$-wtip in $x^{\rho}$.  But, this means that it passes
through two wtips embraced by $x^{\rho}$ of 
ranks less than $\rho$.  But, that in turn means that $S_{1}^{\rho-1}$ and 
$S_{2}^{\rho-1}$ cannot be different $(\rho-1)$-wsections.  $\Box$

An immediate consequence of Lemma 4.1 is 

{\bf Lemma 4.2.}  {\em Any two $\rho$-wnodes $x^{\rho}$ and $y^{\rho}$
that are $\rho$-wconnected but not $\rho$-wadjacent must satisfy
$d(x^{\rho},y^{\rho})\geq \omega^{\rho}$.}

{\bf Theorem 4.3.}  {\em Let
the $\rho$-wsection $S^{\rho}$ of $G^{\nu}$ be locally
$\rho$-finite and have infinitely many boundary $\rho$-wnodes.
Then, given any $\rho$-wnode $x_{0}^{\rho}$ in $S^{\rho}$, there
is a one-ended $\rho$-walk $W^{\rho}$ starting at $x_{0}^{\rho}$:
\[ W^{\rho}\;=\;\langle x_{0}^{\rho},W_{0}^{\alpha_{0}},x_{1}^{\rho},W_{1}^{\alpha_{1}},\ldots, x_{m}^{\rho},W_{m}^{\alpha_{m}},\ldots\rangle \]
such that there is a subsequence of $\rho$-wnodes 
$x_{m_{k}}^{\rho}$, $k=1,2,3,\ldots$, satisfying 
$d(x_{0}^{\rho},x_{m_{k}}^{\rho})\geq\omega^{\rho}\cdot k$.}

The proof of this theorem is just like that of Theorem 10.3 in \cite{gal}
except that the rank 1 therein is replaced by the rank $\rho$ herein. 
In the same way, Corollary 10.4 of \cite{gal} 
generalizes into the following assertion.

{\bf Corollary 4.4}  {\em Under the hypothesis of Theorem 4.3, the 
enlargement $^{*}\!S^{\rho}$ of $S^{\rho}$ has at least one 
$\rho$-hypernode not in its principal galaxy
$\Gamma_{0}^{\rho}(S^{\rho})$ and thus has at least one 
$\rho$-galaxy $\Gamma^{\rho}$ different from its principal $\rho$-galaxy
$\Gamma_{0}^{\rho}(S^{\rho})$.}

\section{When the Enlargement $^{*}\!S^{\rho}$ of a $\rho$-Wsection 
Has a $\rho$-Hypernode Not in the Principal $\rho$-Galaxy of 
$^{*}\!S^{\rho}$}

As always, we take $G^{\nu}$ to be a $\nu$-wgaph with $2\leq\nu\leq
\vec{\omega}$, possibly $\nu=\vec{\omega}$.  
We continue to asssume that the rank $\rho$ 
of a $\rho$-wsection $S^{\rho}$ of $G^{\nu}$ is either a natural
number or $\omega$, but not $\vec{\omega}$.
Let $\Gamma^{\rho}_{a}$ and $\Gamma^{\rho}_{b}$ be two $\rho$-galaxies 
in the enlargement $^{*}\!S^{\rho}$ of a $\rho$-wsection $S^{\rho}$
that are 
different from the principal $\rho$-galaxy $\Gamma^{\rho}_{0}$ of 
$^{*}\!S^{\rho}$.
We shall say that $\Gamma^{\rho}_{a}$ {\em is closer to $\Gamma^{\rho}_{0}$
than is} $\Gamma^{\rho}_{b}$ and that $\Gamma^{\rho}_{b}$ 
{\em is further away from 
$\Gamma^{\rho}_{0}$ than is} $\Gamma^{\rho}_{a}$
if there are a hypernode ${\bf y}=[y_{n}]$ in $\Gamma^{\rho}_{a}$ and a 
hypernode ${\bf z}=[z_{n}]$ in $\Gamma^{\rho}_{b}$ such that, 
for some ${\bf x}=[x_{n}]$ in $\Gamma^{\rho}_{0}$
and for every $m\in\N$, we have
\[ N_{0}(m_{0})\;=\;\{n\!:d(z_{n},x_{n})-d(y_{n},x_{n})
\,\geq\, \omega^{\rho}\cdot m\}\;\in\;{\cal F}. \]
(The ranks of ${\bf x}$, ${\bf y}$, and ${\bf z}$ may have any 
values no larger than $\nu$ other than $\vec{\omega}$.)
Any set of $\rho$-galaxies for which every two of them, say,
$\Gamma^{\rho}_{a}$ and $\Gamma^{\rho}_{b}$  satisfy 
this condition will be said to be
{\em totally ordered according to their closeness to}
$\Gamma^{\rho}_{0}$.  That the conditions for a total ordering 
(reflexivity, antisymmetry, transitivity, and 
connectedness) are fulfilled are readily shown.  
For instance, the proof of Theorem 5.2
below establishes transitivity.

These definitions are independent of the 
representative sequences $\langle x_{n}\rangle$, $\langle y_{n}\rangle$, 
and $\langle z_{n}\rangle$ chosen for ${\bf x}$, ${\bf y}$,
and ${\bf z}$;  the proof of this is exactly the same as 
the proof of Lemma 4.1 of \cite{gal}. 

We will say that a set $A$ is a {\em totally ordered, two-way infinite 
sequence} if there is a bijection from the set $\Z$ of integers to the set 
$A$ that preserves the total ordering of $\Z$.

{\bf Theorem 5.1.}  {\em If the enlargement $^{*}\!S^{\rho}$
of a $\rho$-wsection $S^{\rho}$ of $G^{\nu}$ 
has a hypernode ${\bf v}=[v_{n}]$ (of any rank $\alpha$ with 
$0\leq\alpha\leq\rho$) that is not in 
the principal $\rho$-galaxy $\Gamma^{\rho}_{0}$ of $^{*}\!S^{\rho}$,
then there exists a two-way 
infinite sequence of $\rho$-galaxies totally ordered according to their 
closeness to $\Gamma^{\rho}_{0}$ and with $\bf v$ being in one
of those galaxies.}

{\bf Note.} Here, too, the proof of this is 
much like that of Theorem 4.2 of \cite{gal}, but, since this is the main 
result of this work, let us present a detailed argument.
As always, we choose and fix upon a free ultrafilter $\cal F$.
(Let us also mention that 
a somewhat different version of this theorem with a rather longer 
proof can be found in the archival website, www.arxiv.gov, 
under Mathematics, Zemanian.)

{\bf Proof.}  In this proof, we use the fact that between any two nodes
in a $\nu$-graph there exists a geodesic walk terminating at those nodes;
that is, the length of the walk is equal to the wdistance 
between those nodes.  This is a consequence of the facts 
that the walks terminating at those nodes have ordinal lengths
and that any set of ordinals is well-ordered and thus has a least ordinal.
That least ordinal must be the length of at least one walk terminating 
at those nodes, for otherwise the minimum of the 
walk-lengths would be larger.

As before, let ${\bf x}=[\langle x,x,x,\ldots\rangle]$
be a standard hypernode in $\Gamma_{0}^{\rho}$.
(${\bf x}$ can be of any rank from 0 to $\rho$.)  Since $\bf v$ is not 
in $\Gamma_{0}^{\rho}$, we have that for each $m\in\N$
\begin{equation}
\{n\!: d(x,v_{n})\,>\, \omega^{\rho}\cdot m\}\;\in\;{\cal F}  \label{11.1}
\end{equation}
For each $n\in\N$, if $d(x,v_{n})< \omega^{\rho}\cdot 6$, set $u_{n}=x$, but,
if $d(x,v_{n})\geq \omega^{\rho}\cdot 6$, choose $u_{n}$ such that 
\begin{equation}
d(x,v_{n})\;\leq\; d(x,u_{n})\cdot 3\;\leq\; d(x,v_{n})\cdot 2  \label{11.2}
\end{equation}
That the latter can be done can be seen as follows.

Choose a geodesic $\rho$-walk $W^{\rho}$ terminating at $x$ and $v_{n}$.
Remember that $W^{\rho}$ 
is incident to each of its nonterminal
$\rho$-nodes through at least one $\omega^{\rho-1}$-wtip.
(See Sections 5.3 and 5.5 in \cite{gn} in this regard.)
Moreover, the transition through each $\omega^{\rho-1}$-wtip
contributes $\omega^{\rho}$
to the length of $W^{\rho}$.  Upon tracing $W^{\rho}$ from $x$
toward $v_{n}$, we must encounter at least two $\rho$-nodes, 
both of which are neither closer to $x$ by one-third of the number
of $\omega^{\rho-1}$-wtips traversed by $W^{\rho}$ nor 
further away from $x$ by two-thirds of the number of 
$\omega^{\rho-1}$-wtips 
traversed by $W^{\rho}$. A node on $W^{\rho}$ 
between those two $\rho$-nodes can be
chosen as $u_{n}$.

Suppose there is a $k\in\N$ such that 
$\{n\!: d(x,u_{n})\leq \omega^{\rho}\cdot k\}
\in {\cal F}$. By the left-hand inequality of (\ref{11.2}), 
\[ \{n\!: d(x,v_{n})\,\leq\, (\omega^{\rho}\cdot k)\cdot 3\}\;\supseteq\;
\{n\!: d(x,u_{n})\,\leq\, \omega^{\rho}\cdot k\}\;\in\;{\cal F}. \]
Hence, the left-hand set is a member of $\cal F$, in 
contradiction to (\ref{11.1}).  (These sets cannot both be 
in the ultrafilter $\cal F$.)  Therefore, ${\bf u}=[u_{n}]$
satisfies (\ref{11.1}) for every $m\in\N$ when $v_{n}$ 
is replaced by $u_{n}$;  that is, $\bf u$ is in a galaxy 
different
from the principal $\rho$-galaxy $\Gamma_{0}^{\rho}$.

Furthermore, by the right-hand inequality of (\ref{11.2}),
\[ d(x,u_{n})\;\leq\; (d(x,v_{n})\,-\,d(x,u_{n}))\cdot 2. \]
Suppose there exists a $j\in\N$ such that 
\[ \{n\!: d(x,v_{n})\,-\, d(x,u_{n})\,\leq\, \omega^{\rho}\cdot j\}\;\in\;{\cal F}. \]
Then,
\[ \{n\!: d(x,u_{n})\,\leq\,(\omega^{\rho}\cdot j)\cdot 2\}\;\supseteq\;
\{n\!: d(x,v_{n})\,-\, d(x,u_{n})\,\leq\, \omega^{\rho}\cdot j\}\;\in\;{\cal F} \]
So, the left-hand set is in $\cal F$, in contradiction to our previous 
conclusion that $\bf u$ satisfies (\ref{11.1}) 
with $v_{n}$ replaced by $u_{n}$.  We can conclude that 
$\bf u$ and $\bf v$ are in different $\rho$-galaxies $\Gamma_{a}^{\rho}$
and $\Gamma_{b}^{\rho}$ respectively, with $\Gamma_{a}^{\rho}$ closer
to $\Gamma_{0}^{\rho}$ than is $\Gamma_{b}^{\rho}$.  

We can now repeat this argument with $\Gamma^{\rho}_{b}$ replaced
by $\Gamma^{\rho}_{a}$ to find still another 
$\rho$-galaxy $\tilde{\Gamma}^{\rho}_{a}$ of $^{*}\!S^{\rho}$ different from 
$\Gamma^{\rho}_{0}$ and closer to $\Gamma^{\rho}_{0}$ 
than is $\Gamma^{\rho}_{a}$.  
Continual repetitions yield an infinite sequence of $\rho$-galaxies 
indexed by, say, the negative integers and totally ordered
by their closeness to $\Gamma^{\rho}_{0}$.

The conclusion that there is an infinite sequence of $\rho$-galaxies
progressively further away from $\Gamma^{\rho}_{0}$ than is 
$\Gamma^{\rho}_{b}$ is easier to prove.  
With ${\bf v}\in \Gamma^{\rho}_{b}$ as before,
we have (\ref{11.1}) again for every $m\in\N$. 
Therefore, for each $n\in\N$, we can choose $w_{n}$ as an element of 
$\langle v_{n}\rangle$ such that 
\begin{equation}
d(x,w_{n})\geq d(x,v_{n})
+ \omega^{\rho}\cdot n  \label{4.2}
\end{equation} 
Hence, for each $m\in\N$,
\[ \{n\!: d(x,w_{n})>\omega^{\rho}\cdot m\}\;\supseteq\;\{n\!: d(x,v_{n})>\omega^{\rho}\cdot m\}. \]
Since the right-hand side is a member of $\cal F$, so too is the left-hand side.
Thus, ${\bf w}=[w_{n}]$ is in a galaxy different from $\Gamma_{0}^{\rho}$.
Moreover, from (\ref{4.2}) we have that, for each $m\in\N$,
$\{n\!: d(x,w_{n})-d(x,v_{n})>\omega^{\rho}\cdot m\}$
is a cofinite set and therefore is a member of $\cal F$.
Consequently, $\bf w$ is in a galaxy $\Gamma_{c}^{\rho}$ 
that is further away from $\Gamma_{0}^{\rho}$ than is 
$\Gamma_{b}^{\rho}$.

We can repeat the argument of the last paragraph with
$\Gamma^{\rho}_{c}$ in place of $\Gamma^{\rho}_{b}$ to find 
still another $\rho$-galaxy
$\tilde{\Gamma}^{\rho}_{c}$ further away 
from $\Gamma^{\rho}_{0}$ than is $\Gamma^{\rho}_{c}$.
Repetitions of this argument show that there is an 
infinite sequence of $\rho$-galaxies indexed by, say, the 
positive integers and totally 
ordered by their closeness to $\Gamma^{\rho}_{0}$.  
The conjunction of the two infinite sequences 
along with $\Gamma_{b}^{\rho}$ yields the conclusion of 
the theorem.  $\Box$

By virtue of Corollary 4.4, the 
conclusion of Theorem 5.1 holds whenever $G$ is locally finite.

In general, the hypothesis of Theorem 5.1 may or may not hold.  Thus, 
$^{*}\!S^{\rho}$ either has exactly one $\rho$-galaxy, 
its principal one $\Gamma^{\rho}_{0}$, 
or has infinitely many $\rho$-galaxies.

Instead of the idea of ``totally ordered according to closeness 
to $\Gamma^{\rho}_{0}$,'' we can define the idea of ``partially ordered 
according to closeness to $\Gamma^{\rho}_{0}$'' 
in much the same way.  Just drop the 
connectedness axiom for a total ordering.

{\bf Theorem 5.2.}  {\em Under the hypothesis of Theorem 5.1, 
the set of $\rho$-galaxies of $^{*}\!S^{\rho}$ is 
partially ordered according to the closeness of the $\rho$-galaxies 
to the principal $\rho$-galaxy $\Gamma^{\rho}_{0}$.}

{\bf Proof.}  Reflexivity and antisymmetry are obvious.  
Consider transitivity:  Let $\Gamma^{\rho}_{a}$, 
$\Gamma^{\rho}_{b}$, and $\Gamma^{\rho}_{c}$ 
be $\rho$-galaxies different from $\Gamma^{\rho}_{0}$. (The case where 
$\Gamma^{\rho}_{a}=\Gamma^{\rho}_{0}$ can be argued similarly.)
Assume that $\Gamma^{\rho}_{a}$ is closer to $\Gamma^{\rho}_{0}$ than is 
$\Gamma^{\rho}_{b}$ and that $\Gamma^{\rho}_{b}$ is closer to $\Gamma^{\rho}_{0}$
than is $\Gamma^{\rho}_{c}$.  Thus, for any ${\bf x}$ in $\Gamma^{\rho}_{0}$,
${\bf u}$ in $\Gamma^{\rho}_{a}$, ${\bf v}$ in $\Gamma^{\rho}_{b}$, and 
${\bf w}$ in $\Gamma^{\rho}_{c}$ and for each $m\in\N$, we have
\[ N_{uv}\;=\;\{n\!: d(v_{n},x_{n})-d(u_{n},x_{n})\geq \omega^{\rho}\cdot m\}\,\in\,{\cal F} \]
and 
\[ N_{vw}\;=\;\{n\!: d(w_{n},x_{n})-d(v_{n},x_{n})\geq \omega^{\rho}\cdot m\}\,\in\,{\cal F}. \]
We also have
\[ d(w_{n},x_{n})-d(u_{n},x_{n})\;=\;d(w_{n},x_{n})-d(v_{n},x_{n})
+d(v_{n},x_{n})-d(u_{n},x_{n}). \]
So, 
\[ N_{uw}\;=\;\{n\!: d(w_{n},x_{n})-d(u_{n},x_{n})\geq \omega^{\rho}\cdot 2m\}\,
\supseteq\,N_{uv}\cap N_{vw}\;\in\;{\cal F}. \]
Thus, $N_{uw}\in {\cal F}$. Since $m$ can be chosen arbitrarily,
we can conclude that $\Gamma^{\rho}_{a}$ is closer to $\Gamma^{\rho}_{0}$ 
than is $\Gamma^{\rho}_{c}$.  $\Box$

\end{document}